\documentclass [12pt,reqno]{amsart}
\usepackage{exscale}
\usepackage{graphicx}
\usepackage[centertags]{amsmath}
\usepackage{amsfonts}
\usepackage{ulem}

\renewcommand{\Re}{\operatorname{Re}}
\newcommand{\PGL}{\operatorname{PGL}}
\newcommand{\Hom}{\operatorname{Hom}}

 \newcommand{\sign}{\operatorname{sign}}

\newtheorem*{prop*} {Proposition}
\newtheorem*{conj*} {Conjecture}
\newtheorem*{fact*}{Fact}
\newtheorem*{theorem*} {Theorem}
\newtheorem*{theoremA*} {Theorem A}
\newtheorem*{theoremB*} {Theorem B}

\theoremstyle{remark}

\newtheorem*{defn*}{Definition}
 \newtheorem*{rem*} {Remark}
  \newtheorem*{rems*} {Remarks}

\newtheorem*{acknowledgment*} {Acknowledgment}

\newtheorem*{cl-ex*} {Claim-Example}

\newtheorem*{cl*} {Claim}

\newtheorem*{ex*} {Example}

\newtheorem*{exs*} {Examples}

\newtheorem*{cl-br*} {Claim (B-R)}

\renewcommand{\rho}{\varrho}
\renewcommand{\phi}{\varphi}
\renewcommand{\epsilon}{\varepsilon}

\renewcommand{\rho}{\varrho}
\renewcommand{\phi}{\varphi}

\newcommand{\supp}{\operatorname{supp}}

\newcommand{\GL}{\operatorname{GL}}

\newcommand{\meas}{\operatorname{\meas}}

\newcommand{\nc}{\newcommand}
\nc{\cal}{\mathcal} 
\nc{\la}{\langle} \nc{\ra}{\rangle}
 \nc{\CA}{\cal A}
 \nc{\CBB}{\cal B}
\nc{\CDD}{\cal D}
\nc{\CE}{\cal E}
\nc{\CF}{\cal F} \nc{\CG}{\cal
G} \nc{\CH}{\cal H} \nc{\CI}{\cal I} \nc{\CJ}{\cal J}
\nc{\CK}{\cal K} \nc{\CL}{\cal L} \nc{\CM}{\cal M} \nc{\CN}{\cal
N} \nc{\CO}{\cal O} \nc{\CP}{\cal P} \nc{\CQ}{\cal Q}
\nc{\CR}{\cal R} \nc{\CS}{\cal S} \nc{\CT}{\cal T} \nc{\CU}{\cal
U} \nc{\CV}{\cal V} \nc{\CW}{\cal W} \nc{\CZ}{\cal Z}

\nc{\Ck}{\textsl{k}}

\nc{\fg}{\mathfrak g} \nc{\fii}{\mathfrak i}\nc{\fk}{\mathfrak k}
\nc{\fh}{\mathfrak h} \nc{\fm}{\mathfrak m} \nc{\fn}{\mathfrak n}
\nc{\fA}{\mathfrak A} \nc{\fC}{\mathfrak C} \nc{\fI}{\mathfrak I}
\nc{\fL}{\mathfrak L} \nc{\fS}{\mathfrak S}
\nc{\fz}{\mathfrak z} \nc{\fl}{\mathfrak l}
\nc{\fp}{\mathfrak p}
\nc{\ft}{\mathfrak t}

\nc{\nen}{\newenvironment} \nc{\ol}{\overline}
\nc{\ul}{\underline} \nc{\lra}{\longrightarrow}
\nc{\lla}{\longleftarrow} \nc{\Lra}{\Longrightarrow}
\nc{\Lla}{\Longleftarrow} \nc{\Llra}{\Longleftrightarrow}
\nc{\hra}{\hookrightarrow} \nc{\iso}{\overset{\sim}{\lra}}

\makeatletter
\@addtoreset{equation}{section}
\makeatother
\numberwithin{equation}{section}

 \nc{\ba}{\mathbb A}
 \nc{\bq}{\mathbb Q}
 \nc{\br}{\mathbb R}
 \nc{\bz}{\mathbb Z}
 \nc{\bc}{\mathbb C}
 \nc{\bn}{\mathbb N}
\nc{\bg}{\mathbb G}
 \nc{\ck}{\mathcal{K}}
 \nc{\G}{\Gamma}
 \nc{\sm}{\setminus}
 \nc{\sub}{\subset}
 \nc{\lm}{\lambda}
  \nc{\Lm}{\Lambda}
 \nc{\al}{\alpha}
 \nc{\bt}{\beta}
 \nc{\om}{\omega}
 \nc{\dl}{\delta}
 \nc{\g}{\gamma}
 \nc{\Dl}{\Delta}
 \nc{\Om}{\Omega}
 \nc{\s}{\sigma}
 \nc{\ro}{\rho}
 \nc{\te}{\theta}
 \nc{\SLR}{\operatorname{SL}_2(\br)}
 \nc{\GLR}{\operatorname{GL}_2(\br)}
 \nc{\PGLR}{\operatorname{PGL}_2(\br)}
 \nc{\PSLR}{\operatorname{PSL}_2(\br)}
 \nc{\PSLZ}{\operatorname{PSL}_2(\bz)}
\nc{\PGLZ}{\operatorname{PGL}_2(\bz)}
 \nc{\SLC}{\operatorname{SL}_2(\bc)}
 \nc{\uH}{\mathbb H}
 \nc{\fD}{\mathcal{D}}
 \nc{\fE}{\mathcal{E}}
 \nc{\fO}{\mathcal{O}}
 \nc{\haf}{\frac{1}{2}}
 \nc{\qtr}{\frac{1}{4}}
 \nc{\shaf}{{\scriptstyle\frac{1}{2}}}
 \nc{\hlm}{{\scriptstyle\frac{\lambda}{2}}}

 \nc{\8}{\infty}
 \nc{\7}{{-\infty}}
 \nc{\inv}{^{-1}}
 \nc{\eps}{\varepsilon}
 \nc{\aG}{\mathbf{G}}
 \nc{\spn}{\operatorname{Span}}
 \nc{\Cm}{\operatorname{CM}}
 \nc{\tildl}{\dl^1}
\nc{\chiv}{{\chi_\fp}}
\nc{\psiv}{{\psi_\fp}}
\nc{\piv}{{\pi_\fp}}
\nc{\zt}{Z\setminus T}

\begin{document}

\thispagestyle{plain}

\title{Limiting cycles and periods of Maass forms}

 \author{Andre Reznikov}
 \address{Bar-Ilan University,
Department of Mathematics, Ramat-Gan 52900, Israel}
 \email{reznikov@math.biu.ac.il}
\begin{abstract}We consider (generalized) periods of Maass forms along non-closed geodesics having a closed geodesic as the limit set.
\end{abstract}
\baselineskip18pt
\dedicatory{\it To Leonid Polterovich with best wishes.}
\thanks{The research  was partially supported  by the ERC grant  291612  and by the ISF grant 533/14.}
\maketitle
\section{Introduction}

The aim of this note is to introduce a new kind of a period for Maass forms. For outsiders of the automorphic world, we recall that in the theory of automorphic functions eigenfunctions of the Laplacian on a finite volume hyperbolic Riemann surface are called Maass forms (after H. Maass who realized their importance in Number Theory). There are 3 types of closed cycles naturally appearing in the theory of automorphic functions on the group $G=\PGL_2(\mathbb R)$ with respect to a lattice $\G\subset G$. These are closed horocycles which are associated to closed orbits in the automorphic space $X=\G\sm G$ of a unipotent subgroup $N\subset G$, closed geodesics and geodesic rays starting and ending in a cusp (both of these types are associated to closed orbits on $X$ of the diagonal subgroup $A$), and closed geodesics circles which are associated to orbits on $X$ of a maximal compact subgroup $K$.

We propose to consider one more period along a special type of a {\it non-closed} orbit of the subgroup $A,$ i.e., along a {\it non-closed} special geodesic on the corresponding Riemann surface. These geodesics will have closed geodesics as their limit sets. Our justification for introducing such cycles is that (generalized) periods of Maass forms along these geodesics satisfy nice analytic properties.

Periods (with characters) along ``classical'' cycles lead to Fourier coefficients of cusp forms and to $L$-functions (e.g., the Hecke $L$-function given by a period along the geodesic ray connecting two cusps and a special value of a quadratic base change $L$-function from a theorem of Waldspurger \cite{Wa} appearing as a period along a closed geodesic), and hence play an important role in Number Theory. Admittedly, we do not know yet what is the arithmetic meaning of these new periods (although their residues are connected to periods along closed geodesics and hence to special values of $L$-functions).

\subsection{Limiting cycles} Although our methods do {\it not} use arithmetic, we describe the phenomenon in the simplest (and probably in the most interesting) case of the modular lattice $\Gamma =\PGL_2(\mathbb Z)\subset G.$ Denote  by $Y=X/K=\G\sm\CH$  the corresponding Riemann surface (here $K=O(2)$ and $\mathcal H=G/K$ denotes the hyperbolic upper half plane).  Let $\ell\subset Y$ be a closed geodesic. In particular, lifting $\ell$ to the upper half plane $\mathcal H$, we obtain a geodesic $\tilde\ell\subset \mathcal H$ connecting points $\alpha,$ $\overline\alpha$ on the absolute $\mathbb P^1(\br)$ of $\mathcal H$ (for the lattice $\PGL_2(\mathbb Z)$, $\alpha$ and $\overline \alpha$ are two conjugate real quadratic irrationalities).

\vskip 8cm
\includegraphics[scale=0.15]{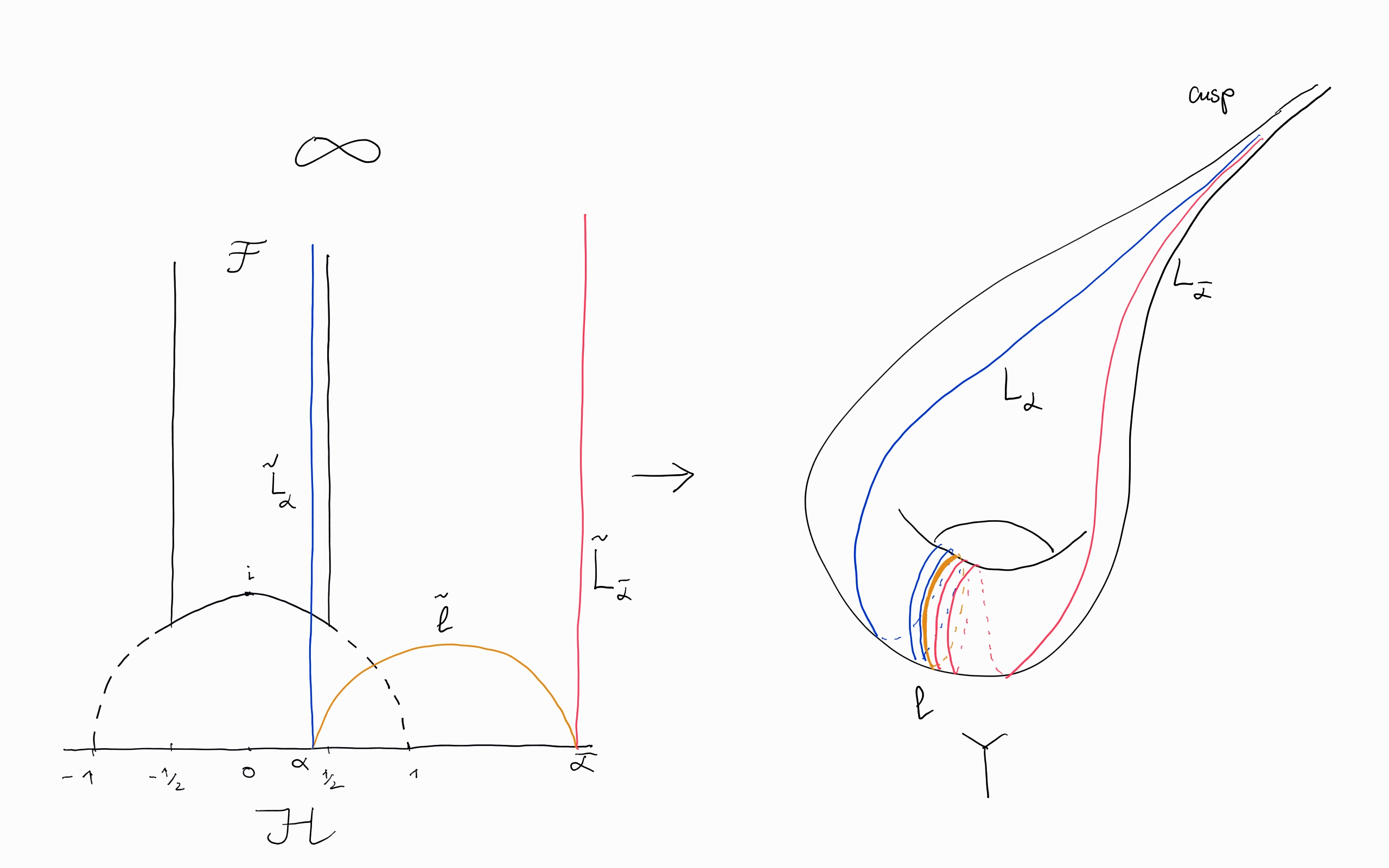}

We consider the vertical geodesic $\tilde L=\tilde L_\al\subset\mathcal H$ connecting the cusp of $\G$ at $\infty,$ and let say $\alpha$. On $Y$ the corresponding geodesic $L\subset Y$ becomes a geodesic winding around the closed geodesic $\ell$ in one direction and escaping into the cusp in the opposite direction. Clearly, $\ell$ is in the closure of $L$ (in fact, $\overline L=L\cup\ell).$ In general, $\ell$ is a piecewise smooth geodesic in the case it passes through  a conical point on $Y$ (i.e., through a fixed point of an elliptic element in $\G$).

\subsection{Limiting periods of Maass forms} Let $\varphi=\varphi_\lambda$ be a cuspidal Maass form on $Y$ with the eigenvalue $\mu=\frac{1-\lambda^2}{4}$ of the Laplacian $\Dl$ on $Y$, where $\lambda\in\mathbb R\cup(0,1)$.
For $s\in\bc$, we consider a twisted period of $\varphi$ along $L,$ namely,
\begin{equation}\label{eq1}
I_s(\varphi)=I_{s,L_\alpha}(\varphi_\lm)=\int_L\varphi(l)|l|^sdl=\int_{\tilde L}\varphi(\tilde l)|\tilde l|^sd\tilde l,
\end{equation}
where we identify $L$ with the subgroup $A_{\alpha,\infty}=\text{Stab}_G(\tilde L)$ (e.g., by choosing a point on $\tilde L$), and, using this identification, consider multiplicative characters $|\cdot|^s$ on $\tilde L \simeq  \mathbb R_+^\times,$ and an $A_{\alpha,\infty}$ invariant measure $d\tilde l$ on $\tilde L$. We can write the same integral in ``additive'' coordinates on $L$ using the natural geodesic parameter $t\in\br$ on $L$ (i.e., $t=\ln a$) and the corresponding measure $dt$ given by the length.  We have then $I_s(\varphi)=\int_L\varphi(t)e^{st}dt$.

It is easy to see that the integral $I_s(\varphi)$ is absolutely convergent for $\Re(s)>1$ (e.g., using the fact that $\phi$ is a bounded function with an exponential decay in  cusps).

Our main result is the following

\begin{theoremA*}  The integral $I_s(\varphi)$ is absolutely convergent for $\Re(s)>1$, and has meromorphic continuation to $\mathbb C$.
\end{theoremA*}

Moreover, we will identify possible poles and residues of $I_s(\varphi)$ in terms of periods of $\varphi$ along the closed geodesic $\ell\subset Y.$ We also note that  non-cuspidal Maass forms could be also treated by our method.

One can give a reformulation of Theorem A  in terms of Fourier coefficients of $\varphi.$ Let \begin{equation}\label{eq2}
\varphi_\lambda(x+iy)=\sum_{n\ne0} a_n\cdot K_{\lambda,n}(Y)e^{2\pi inx}
\end{equation}
be the Fourier expansion of $\varphi_\lambda$ at  $\infty.$ Here we denote by $K_{\lambda,n}(Y)$ a variant of the classical $K$-Bessel function which we will normalize by the correspong matrix coefficient (see \eqref{Bessel}). Coefficients $a_n=a_n(\varphi_\lambda),$ up to a simple twist, coincide with classical Fourier coefficients of $\varphi_\lambda$, and, in particular,  satisfy the asymptotic $\sum\limits_{|n|\leq T}|a_n|^2 \sim C\cdot T$ with some constant $C=C(\varphi_\lambda)>0$ (of course, this depends on normalization of functions $K_{\lm,n}).$

By writing the function $\phi$ on $\tilde L$ through its Fourier expansion at the cusp, we obtain a Dirichlet series expression for the integral $I_s(\phi)$. This leads to the  following, more ``classical'', form of Theorem A. 
\begin{theoremB*}  The twisted Dirichlet series
\begin{equation}\label{eq3}
D_{s,\alpha}(\varphi)=\sum_{n\ne0} \frac{a_n}{|n|^s}\cdot e^{2\pi i n\alpha }
\end{equation}
is absolutely convergent for $\Re(s)>1,$ and has meromorphic continuation to $\mathbb C.$
\end{theoremB*}

\begin{rems*}

1. Comparing poles and residues of $D_{s,\alpha}$ (which are given in terms of twisted periods of $\varphi$ along $\ell$), we see that the combination $D_{s,\alpha}-D_{s,\overline \alpha}$ is holomorphic (or what is the same, the combination $I_{s,L_\alpha}(\varphi_\lambda)+I_{s,\overline L_{\overline\alpha}}(\varphi_\lambda)$
is holomorphic where $\overline L_{\overline \alpha}$ is the cycle connecting $\overline \alpha$ to the cusp at $\infty$ and ``$\overline{L}$'' stands for the cycle with the reversed orientation). 

Our method also naturally  treats the signed sum \begin{equation}\label{eq3-}
D^\pm_{s,\alpha}(\varphi)=\sum_{n\ne0} \frac{a_n}{|n|^s}\cdot e^{2\pi i n\alpha }\sign(n)\ ,
\end{equation} and as a result the sum  $D^>_{s,\alpha}(\varphi)=\sum_{n>0} \frac{a_n}{|n|^s}\cdot e^{2\pi i n\alpha }$.

Our proof also gives a polynomial in $s$ and $\lm$ bound for $D_{s,\alpha}(\varphi_\lm)$. 

2.  Dirichlet series for the ``additively twisted'' $L$-function
$$D_{s,a}(\varphi)=\sum_{n\ge1}\frac{a_n(\phi)}{n^s}e^{2\pi ina},\quad a\in\mathbb R$$
were studied extensively starting with Hecke for $a\in\mathbb Q$ (this leads to the usual Hecke $L$-function and hence have holomorphic  continuation to $\bc$). Properties of $D_{s,a}$ for an irrational $a$ were studied, among other places, in \cite{C} and \cite{MS}, and in \cite{M}   for $GL(n)$, and it is certainly  expected that for a general $a$, the corresponding series does not have analytic continuation to $\bc$. 

3. The proof we give does not use arithmetic. It is interesting to see if for the arithmetic $\Gamma $ and Hecke-Maass forms, the above Dirichlet series have some arithmetical meaning (from what we have seen above, residues of this series do have a connection to special values of $L$-functions via the theorem of Waldspurger \cite{Wa}).

We note that the cycle over which we integrate has clear arithmetic description (for a congruence subgroup of $\G=\PGLZ$). Let $E$ be the quadratic field corresponding to the geodesic $\ell$ (in particular,  the primitive element $\g_\ell$ corresponds to a unit in an order in $E$). We have the natural imbedding $E^\times\hookrightarrow\GL_2(\bq)$ which gives rise to the real tori $T_E(\br)=E^\times(\br)\hookrightarrow\GL_2(\br)$. In the quotient space this gives the geodesic $\ell$ (or rather a collection of geodesics corresponding to the class group of the corresponding order). The corresponding tori $T_E(\br)$ has $\al$ and $\bar\al$ as fixed points on $\mathbb P^1(\mathbb R)$. In the light of the Waldspurger theorem \cite{Wa} which deals with the period squared $|\int_\ell\phi|^2$, it is natural to consider the imbedding $T_E(\br)\times T_E(\br)\hra \GL_2(\br)\times\GL_2(\br)$ coming from two conjugated imbeddings $E\hra\br$. Note that pointwise these two imbeddings of $T_E(\br)$ coincide with one fixing the pair $\al,\bar\al$ and another the pair $\bar\al,\al$. In the quotient space $X\times X$ the corresponding cycle $\ell\times\ell$ is compact. Our period $L_\al$ or rather the product $L_\al\times L_{\bar\al}\subset X\times X $ could be described in similar terms. Consider another tori $T'_E\subset \GL_2(E)$ which fixes $\8$ and $\al$ on $\mathbb P^1(E)$. $T'_E$ is defined over $E$ while $T_E$ is defined over $\bq$ (both are conjugate in $\GL_2(E)$). $T'_E$  has two conjugated imbeddings into $\GL_2(\br)$. Hence we obtain the imbedding $T'_E(\br)\times T'_E(\br)\hra  \GL_2(\br)\times \GL_2(\br)$ with one component fixing the pair $\8,\al$ and another $\8,\bar\al$. In the quotient space this becomes our cycle $L_\al\times L_{\bar\al}$. The corresponding cycle is not compact since there are no points in $T'_E$ which are defined over $\bq$ (unlike for the tori $T_E$). This shows that the product of integrals $I_s(\phi)\cdot I_{-s}(\bar\phi)$ is given by an integral over a  cycle in $Y\times Y$ which is defined over $E$ while the Waldspurger cycle is defined over $\bq$.

4. For another instance of an analytic continuation of a pairing $\langle I,v_s\rangle$ for a singular vector $v_s\in V_\lm$, see \cite{R2}. In that case, singularities of the test vector $v_s$ are at rational points on $\mathbb P^1(\mathbb R)$, and the analytic continuation is achieved with the help of Hecke operators.

5. Integrals along limiting cycles $
L_\alpha$ were considered in \cite{MM} from a different point of view for the group $\PSLZ$. In particular, it is shown in \cite{MM} that for any fixed $x\in \tilde L_\al$ and  $y\in \tilde L_\al$, the mean value of the  integral $I(\phi,\al,x)=\lim\limits_{y\to \al}|\tilde L_{x,y}|\inv\int_{\tilde L_{x,y}}\phi(l)dl$  converges to the period integral $\int_\ell\phi(l)dl$ for $\al$ which is a fixed point of a hyperbolic element in $\PSLZ$ (and in particular is independent of $x$), and vanishes for other irrational $\al$. For $\al$ which is a fixed point of a hyperbolic element in $\PSLZ$, we can deduce this result from the existence of the simple pole of $I_s(\phi)$ at $s=1$ with the residue given by the period along the closed geodesic $\ell$. In fact, we can treat also twisted integrals $\lim\limits_{y\to \al}|\tilde L_{x,y}|\inv\int_{\tilde L_{x,y}}\phi(l)\chi(l)dl$ where, as in \eqref{eq1}, we view $L_\al$ as an orbit of $A$ and $\chi:A\to\bc^\times$ a character trivial on $\G_\ell=A\cap \G$. Looking at poles of  $I_s(\phi)$ at appropriate $s$, we see that the limit is given by the twisted period $\int_\ell\phi(l)\chi(l)dl$ along the corresponding closed geodesic $\ell$. 

6. We also have analogous results for integrals over limiting geodesics winding down from one closed geodesic to another closed geodesic (such a geodesic has two rotational numbers, for a general discussion, see  \cite{P}), although these integrals are not easily expressible in terms of Fourier coefficient (and in fact, such integrals could be treated for co-compact lattices as well).

We  note that there are infinitely many non-homotopic geodesic having two closed geodesic (or a closed geodesic and a cusp) as their limit set. Our result is valid for any such a geodesic. 

\end{rems*}
\subsubsection*{Acknowledgments} It is a pleasure to thank Joseph Bernstein for numerous discussions. 

 \section{Proof}

 We now present the proof of Theorem B (and Theorem A follows easily from Theorem B by writing $\varphi$ on $L_\alpha$ through its Fourier expansion).

 Our proof is based on standard methods in representation theory of automorphic functions (and, in particular, use of co-invariants of representation of $\PGL_2(\mathbb R)$ with respect to some cyclic subgroups coming from $\G$).

 \subsection{Representation theoretic setup}

 Let $\pi_\lambda$ be an automorphic representation associated to a cuspidal Maass form $\varphi_\lambda$ (this means that we have an abstract unitary representation $\pi_\lambda$ of $G$ and an isometry $\nu: L_{\pi_\lambda}\to L^2(X)$ such that $\nu(e_0)=\varphi_\lambda$ for a $K$-fixed vector $e_0\in V_{\pi_\lambda}$ of norm one). We will work with the space of smooth vectors $V_\lambda=L_{\pi_\lambda}^\infty\subseteq L_{\pi_\lambda}$ and will assume that $\lambda\in i\mathbb R$ (i.e., $\pi_\lambda$ is a principal series representation of $G$). Frobenius reciprocity (e.g., see \cite{BR}) gives the corresponding functional $I=I_{\varphi_\lambda}\in\Hom_\lambda(V_\lambda,\mathbb C)$, i.e., a $\Gamma$-invariant functional $I: V_\lm\to\bc$ such that  $\varphi_v(g):=\nu(v)(g)=I(\pi_\lambda(g)v)$ for any smooth vector $v\in V_\lambda.$ The functional $I$ belongs to some Sobolev class of functionals on $V_\lambda.$

 In \cite{BR}, we showed that the functional $I$ has finite $-\left(\frac{1}{2}+\varepsilon\right)$-Sobolev $L^2$-norm for any $\eps>0$ (similar results for the H\"older class were obtained in \cite{O}, \cite{S}). The space $V_\lambda$ has natural realization in the space of smooth functions $C_\lambda^\infty(\mathbb R)$ on $\mathbb R$ with some condition on the decay at infinity (e.g.,  functions have asymptotic $|x|^{-1+\lm}$ as $|x|\to\8$). We call this model the linear model of $\pi_\lambda.$ Hence we can view $I$ as a distribution on $\mathbb R.$ The element $\left(\begin{smallmatrix} 1 & 1\\ 0 & 1\end{smallmatrix}\right)\in\G$ generates the cuspidal subgroup $\G_\8=\G\cap N$, and acts by the usual translation $\pi_\lambda \left(\begin{smallmatrix} 1 & 1\\ 0 & 1\end{smallmatrix}\right)f(x)=f(x-1)$ in the linear model of $V_\lambda$. Hence the $\Gamma$-invariant distribution $I$ has the periodic Fourier expansion
 \begin{equation}\label{eq4}
 I=\sum_{0\ne n\in\mathbb Z} a_nW_{n,\lambda},
 \end{equation}
 where $W_{n,\lambda}:V_\lambda\to{\mathbb C}$ is a Whittaker functional. We normalize Whittaker functionals by their norm on the space of $\Gamma_\infty$-co-invariants $(V_\lambda)_{\Gamma_\infty}\simeq C^\8([0,1])$ in the following way. We view $W_{n,\lambda}$ as a generalized vector in $\pi_{-\lambda}$. We fix $W_{1,\lambda}(x)=e^{2\pi ix}$ and denote by $W_{n,\lambda}(x)=|n|^{-\haf}\cdot\pi_{-\lambda}\left(\begin{smallmatrix} n^{-\haf} & 0\\ 0 & n^{\haf}\end{smallmatrix}\right)W_{1,\lambda}=|n|^{-\frac{\lambda}{2}}\cdot e^{2\pi inx}$. Hence up to a factor of the  absolute norm one, $W_{n,\lambda}$ is given by the exponent $e^{2\pi inx}.$

Note that with such a normalization of functionals $W_{n,\lambda}$, coefficients $a_n$ in \eqref{eq4} do not depend on the change of sign of $\lambda$ (but do depend on a choice of $W_{1,\lambda}).$

In order that coefficients $a_n$ in \eqref{eq4} will be consistent with the Fourier coefficients of $\varphi_\lambda$ in \eqref{eq2}, we normalize $K$-Bessel functions in \eqref{eq2}, by the corresponding matrix coefficient, i.e., by
\begin{equation}\label{Bessel} K_{\lambda,n}(y)=\Big\langle \pi_\lambda\left(\begin{smallmatrix} y^{\haf} & 0\\ 0 & y^{-\haf}\end{smallmatrix}\right)e_0,\, W_{n,\lambda}\Big\rangle=|n|^{-\lambda/2}|y|^{\frac{\lambda-1}{2}}\int_{\mathbb R}|1+(y^{-1}x)^2|^{\frac{\lambda-1}{2}}e^{-2\pi inx}dx\ ,\end{equation}
where $e_0(x)=|1+x^2|^{\frac{\lambda-1}{2}}$ is a $K$-fixed vector in the linear model of $V_\lambda.$

From \eqref{eq4} we see that the Dirichlet series $D_{s,\alpha}$, up to a ratio of $\G$-functions, is given  by the pairing of $I$ with a (generalized) vector $d_{s,\alpha}$  given by the kernel $d_{s,\alpha}(x)=|x-\alpha|^{-1-\frac{\lambda}{2}+s}$ in the linear model of the representation $\pi_\lambda$. 
We have
$$\langle I,d_{s,\alpha}\rangle=\sum_{0\ne n\in\mathbb Z}a_n\cdot W_{n,\lambda}(d_{s,\alpha})=\gamma(s-\shaf\lambda)\cdot \sum_{n\ne0}a_n\cdot n^{-s}\cdot e^{2\pi in\alpha}=\gamma(s-\shaf\lambda)\cdot D_{s,\alpha(\varphi)},$$
where the Fourier transform gives $\mathcal F(|x|^{s-1})[\xi]=\gamma(s)\cdot|\xi|^{-s} $ with $\gamma(s)=\pi^{-\frac{s}{2}}\Gamma\left(\frac{s}{2}\right)\big/ \pi^{-\frac{1-s}{2}} \Gamma
\left(\frac{1-s}{2}\right).$

We remark that it is more natural to shift $s$ by $\frac{1}{2}$ and consider the vector (which we will denote by the same letter) $d_{s,\alpha}(x)=|x-\alpha|^{-\frac{1}{2}-\frac{\lambda}{2}+s}.$

The vector $d_{s,\alpha}$ is not a smooth vector, and we have to make sense out of the pairing $\langle I,d_{s,\alpha}\rangle.$ We will show that in fact this pairing is well-defined for $\Re s>1$ (since the corresponding vector belongs to an appropriate Sobolev space) and, moreover, the corresponding pairing has meromorphic continuation (because $d_{s,\alpha}$ belongs to an appropriate Sobolev space after the projection to co-invariants of $V_\lambda$ with respect to the cyclic subgroup $\Gamma_\ell\subset \G$ corresponding to the geodesic $\ell$).

\subsection{Co-invariants}

Our aim is to analytically continue the expression $\langle I,d_{s,\alpha}\rangle,$ where $d_{s,\alpha}(x)=|x-\alpha|^{-\frac{1}{2}-\frac{\lm}{2}+s}.$ It is easy to see that $d_{s,\alpha}$ is an eigenvector of the action of the torus $A_{\alpha,\infty}\subset G$  on $\mathbb P^1(\mathbb R)$. Namely, we have $\pi_\lambda(a_{\alpha,\infty})d_{s,\alpha}=|a_{\alpha,\infty}|^{-2s}\cdot d_{s,\al}.$
Here $|a_{\alpha,\infty}|^{-2s}$ denotes the natural character $\chi_{-2s}:A_{\alpha,\infty}\to\mathbb C^\times,$ $s\in\mathbb C$ induced by the conjugation of $A_{\alpha,\infty}$ to $A_{0,\infty}=A=\left\{\left(\begin{smallmatrix} a & 0\\ 0 & a^{-1}\end{smallmatrix}\right)\ |\ a\in\br^\times\right\} $ and the corresponding (quasi-) character $\chi_\sigma:A_{0,\infty}\to\mathbb C^\times$ is given by $\chi_\sigma\left(\begin{smallmatrix} a & 0\\ 0 & a^{-1}\end{smallmatrix}\right)=|a|^\sigma$ for $\s\in\bc$.

We now use the fact that $\text{Stab}_\Gamma(\infty)=\Gamma_\infty\ne\{e\}$ and that $\text{Stab}_\Gamma(\alpha)=\Gamma_\ell\ne\{e\}$. (Note that $\text{Stab}_\Gamma((\alpha,\infty))=\Gamma\cap A_{\alpha,\infty}=\{e\}.)$

The generalized vector $d_{\alpha,s}$ is smooth outside of $\infty$ and $\alpha.$ We use partition of unity to separate these singularities. Let $u_0$ be a real valued function supported near $\alpha$ and $u_\infty$ the complimentary function (i.e., $1\ge u(x)\ge0,$ $\supp u\subseteq[\alpha-1,\alpha+1],$ $u\big|_{[\alpha-1/2,\alpha+1/2]}\equiv 1$ and $u_0+u_\infty\equiv1).$
Denote by $d_{\alpha,s}^0(x)=u_0(x)d_{\alpha,s}(x),$ $d_{\alpha,s}^\infty(x)=u_\infty(x)d_{\alpha,s}(x).$

\subsubsection{Claim A} $\langle I,d_{\alpha,s}^\infty\rangle$ is holomorphic.

This follows from the fact that in co-invariants $(V_\lambda)_{\Gamma_\infty}=V_\lambda\big/\{(1-n)f\}_{n\in\Gamma_\infty,f\in V_\lambda},$ the vector $d_{\alpha,S}^\infty$ is differentiable $N$-times for any $N\ge1$ (and, in particular, it belongs to the $\left(\frac{1}{2}+\varepsilon\right)$-Sobolev space of $V_\lambda$ for any $\eps>0$).

\subsubsection{Claim B} $\langle I,d_{\alpha,s}^0\rangle$ is well-defined for $\Re(s)>\frac{1}{2},$ and has the meromorphic continuation to $\mathbb C$.

For $\Re(s)>\frac{1}{2},$ the function $d_{\alpha,s}^0(x)=|x-\alpha|^{-\frac{1}{2}-\frac{\lambda}{2}+s}u_0(x)$ has a polynomial singularity at $x=\alpha$ which belongs to a $\left(\frac{1}{2}+\varepsilon\right)$-Sobolev space for $\Re(s)>\haf +\eps$. Hence the value of $\langle I,d_{\alpha,s}^0\rangle$ is well-defined and holomorphic in $s$ in this domain.

Consider a non-trivial element $\gamma_\ell\in\Gamma_\ell=\Gamma\cap A_{\alpha,\overline\alpha}$ (e.g., a primitive hyperbolic element in $\Gamma$ which corresponds to the geodesic $\ell$). The element $\gamma_\ell$ is conjugated to a diagonal element $a_\ell=\left(\begin{smallmatrix} q_\ell^{\frac{1}{2}} & 0\\ 0 & q_\ell^{-\frac{1}{2}}\end{smallmatrix}\right)\in A$ for some $q_\ell>1$ ($q_\ell$ is easily given in terms of length of $\ell$). We denote by $A_\ell=\langle a_\ell\rangle\subset A$ the corresponding subgroup. The space of co-invariants $(V_\lambda)_{\Gamma_\ell}$ is naturally isomorphic to the space of (smooth) functions on $A_\ell\setminus A\simeq\ \langle q\rangle\ \sm \mathbb R^\times.$ Up to a conjugation (by an  element $g_\ell$ conjugating $\gamma_\ell$ to $a_\ell$), our function $d_{\alpha,s}^0$ becomes the function $\tilde d^0_{\beta,s}(x)=\pi_\lm(g_\ell)d_{\alpha,s}^0=|x|^{-\frac{1}{2}-\frac{\gamma}{2}+s}|x-\beta|^{-\frac{1}{2}-\frac{\lambda}{2}-s}\tilde
u_0(x))$ with $\tilde u_0$ supported in an interval which includes $0,$ but not $\beta$ (here $\beta=g_\ell(\8)$ and $0=g_\ell(\al)$), and $\tilde u_0\equiv 1$ in some interval which includes $0.$ Hence $\tilde d^0_{\beta,s}$ has the only singularity at $x=0.$

To describe the image of $\tilde d^0_{\beta,s}$ in the space of co-invariants $(V_\lambda)_{\langle a_\ell\rangle}\simeq C^\infty(\langle q\rangle\setminus\mathbb R^\times),$ we can compute (multiplicative) Fourier coefficients $b^+_j(s)=\Big\langle \tilde d^0_{\beta,s},|x|^{-\frac{1}{2}-\frac{\lambda}{2}+s_j}\Big\rangle_{V_\lambda\simeq L^2(\mathbb R)}$ where $s_j\in i\mathbb R$ are defined by the condition $\pi_\lambda\left(\begin{smallmatrix} q_\ell{\frac{1}{2}} & 0\\ 0 & q_\ell^{-\frac{1}{2}}\end{smallmatrix}\right)|x|^{-\frac{1}{2}-\frac{\lambda}{2}+s_j}=|x|^{-\frac{1}{2}
-\frac{\lambda}{2}+s_j}$ (namely,  $|x|^{-\frac{1}{2}-\frac{\lambda}{2}+s_j}$ are $A$-equivariant functionals on $V_\lambda$ which are trivial on $A_\ell=\langle a_\ell\rangle$), i.e., $s_j=\frac{2\pi i}{\ln q}\cdot j,$ $j\in\mathbb Z.$

In fact, since $A\simeq\br^\times$ has two connected components, we also have to consider ``odd'' multiplicative Fourier coefficients $b_j^-(s)=\Big\langle \tilde d^0_{\beta,s},|x|^{-\frac{1}{2}-\frac{\lambda}{2}+s_j}\sign(x)\Big\rangle_{V_\lambda\simeq L^2(\mathbb R)}$.

We claim that coefficients $b_j^\pm(s)$ are meromorphic in $s$, and for a fixed $s$ outside of poles, are rapidly decreasing as $|j|\to\8$. This implies that the vector $d_{\alpha,s}^0$ has a smooth representative in co-invariants (outside of poles) and hence could be paired with $I$. 

Taking first $N$-terms of the Taylor expansion of the {\it smooth} function $|x-\beta|^{-\frac{1}{2}-\frac{\lambda}{2}-s}$ at $x=0$, we see that
$$F_{0,\beta,s}(x)=|x-\beta|^{-\frac{1}{2}-\frac{\lambda}{2}-s}=a_0(\lambda,s)+a_1(\lambda,s)x+\dots+O_{s,\lm,\bt}(|x|^N).$$ Hence it is enough to 
compute scalar products in $V_\lambda\simeq L^2(\mathbb R)$ of the form $\Big\langle |x|^{-\frac{1}{2}-\frac{\lambda}{2}+s+n}\tilde u_0(x),\, |x|^{-\frac{1}{2}-\frac{\lambda}{2}+s_j}\sign^\epsilon(x)\Big\rangle,$ $n\ge 0$, $\eps=0,1$. These are well-defined for $\Re(s)\gg0$ and have meromorphic continuation to $\mathbb C$ with poles appearing when
$-\frac{1}{2}-\frac{\lambda}{2}+s+n-\frac{1}{2}+\frac{\lambda}{2}-s_j=k$ is a
{\it negative} integer. 
 The residue at such a point is given by the coefficients $\rho_{s_j,\eps}=\big\langle I,d_{s_j,\Gamma_\ell}\rangle_{(V_\lambda)_{\Gamma_\ell}},$ where $d_{s_j,\eps\G_\ell}$ is the $A_{\alpha,\overline\alpha}\ $-equivariant functional $d_{s_j,\Gamma_\ell}:V_\lambda\to\mathbb C$ satisfying $\pi_{-\lambda}(a)d_{s_j,\eps,\Gamma_\ell}=|a|^{s_j}\sign^\eps(a)d_{s_j,\eps,\Gamma_\ell},$
$a\in A_{\alpha,\overline\alpha}$ with the natural meaning of the character $|a|^s\sign^\eps(a)$ (recall that $A_{\alpha,\overline\alpha}=A_\ell $ is the subgroup which is conjugate to the diagonal subgroup $A=A_{0,\infty}$ and such that $\Gamma\cap A_\ell=\Gamma_\ell$).  Coefficients $\rho_{s_j,\eps}$ are nothing else but (appropriately normalized) twisted (by the character $|a|^{s_j}\sign^\eps(a)$) periods of $\varphi$ along the geodesic $\ell$. In fact, since we deal with the group $\PGLR$ and the corresponding diagonal subgroup $A$ has two connected components, orbits of $A$ on $X$ are coming in pairs of connected orbits or one connected orbit with an automorphism (which could be trivial). Periods 
$\rho_{s_j,\eps}$ reflect this structure (this is discussed in detail in \cite{R1}). 

We note that it is easy to see that coefficients $b_j^\pm(s)$ are polynomial in $s$ (and in $\lm$) and this gives polynomial bounds for the analytically continued series $D_{s,\al}(\phi_\lm)$.

\end{document}